\input amstex
\input xy
\xyoption{all}
\documentstyle{amsppt}
\document
\magnification=1200
\NoBlackBoxes
\nologo
\hoffset=0.4in
\voffset=0.7in
\def\r{\roman}

\def\Z{\bold{Z}^+}
\def\R{\bold{R}}
\def\C{\bold{C}}
\def\c{\Cal{C}}
\def\H{\Cal{H}}
\def\N{\bold{N}}
\vsize16cm


\bigskip

 \centerline{\bf  RENORMALIZATION AND COMPUTATION II:}
 
 \medskip
 
 \centerline{\bf TIME CUT--OFF AND THE HALTING PROBLEM}

\medskip

\centerline{\bf Yuri I. Manin}

\medskip

\centerline{\it Max--Planck--Institut f\"ur Mathematik, Bonn, Germany,}

\centerline{\it and Northwestern University, Evanston, USA}

\bigskip

{\bf Abstract.} This is the second installment to the  project initiated in [Ma3].
In the first Part, I argued that both philosophy and technique
of the perturbative renormalization in quantum field theory could be meaningfully
transplanted to the theory of computation, and sketched several
contexts supporting this view.

In this second part, I address some of the issues raised in [Ma3] and 
provide  their development in three contexts: a categorification of the algorithmic
computations; time cut--off and Anytime Algorithms; and finally, a
Hopf algebra renormalization of the Halting Problem.

\bigskip

\centerline{\bf  Contents}

\medskip

0. Introduction

1. Enriched programming methods with typing and parallelism:
a categorical approach

2. Cut-off regularization and Anytime Algorithms

3. Regularization and renormalization of the Halting Problem

\bigskip

\centerline{\bf 0. Introduction}

\medskip

{\bf 0.1. Regularization and Anytime Algorithms.}
It is well known that classical theory of computability
includes as its organic part phenomena of 
{\it non}--computability. Namely, an attempt to compute the value of
a partially recursive function at a point where it is not defined,
might stall the computation forever,
but we will never know whether this is so or simply we did not
wait long enough  (``the Halting Problem is undecidable'').
\smallskip

Applied theory of computation deals with algorithms
processing finite amount of data into finite outputs.
Nevertheless, even in such
theoretically safe situations time/memory requirements
may make the implementation of a sound algorithm unfeasible. 

\smallskip

The celebrated theory of polynomial time computations
and discovery of the $P/NP$--problem served
as a neutral zone meeting point between theoretical possibility and practical
feasibility, and revealed beautiful new mathematical
structures.

\smallskip

However, applied computer scientists consider other
possible ways of turning unfeasible computations into feasible ones,
known under the code word ``Anytime Algorithms''.
Basically, an ``Anytime Algorithm'' allows the computation
to stop at a feasible time, and supplies the result
of such a mutilated procedure with a measure of its quality.
See [GrZi] and a nice short introduction  [Gr].

\smallskip

In the Sec. 2 of this article, the second installment to the project initiated in [Ma3],
it is suggested that theoretically ``Anytime Algorithms''
can be treated as one of the versions of
regularization schemes
in Quantum Field Theory: {\it time cut--off} (for more detailed description
of the whole project, see Introduction
to [Ma3]). More precisely, I analyze from this viewpoint results
of the stimulating paper  by Ch.~Calude, M.~Stay ``Most programs stop quickly or never halt'' ([CalSt1]).

\smallskip

One of the themes, that the analogies
with renormalization and experience with
Anytime Algorithms bring to the foreground in the computation theory,
is the stress on  the structure of programs determined by the operation
{\it ``composition of programs''} and by the  {\it explicit parallelism}, 
that played a key role
in our treatment of perturbative renormalization as a 
model for regularizing computations in [Ma3].

\smallskip   

Notice that many standard descriptions of
programming methods are not stable with respect to the
composition and have no natural means for expressing
parallelism. 

\smallskip
For example, composition $T_2\circ T_1$ of two Turing machines,
informally defined as computation in which the output of $T_1$
becomes the oracular input (``program'') for $T_2$, is not directly described
as a new Turing machine $T_3$. 

\smallskip
Language--like constructions
such as lambda--calculus, being inherently linear/se\-quential,
are not well--suited for expressing options of parallelism.

\smallskip

"Flowcharts'' imagery for which I made some propaganda 
in [Ma3] serves these goals much better. In the Section 1 of this article,
I show that the same ideas admit a succinct
categorical expression, and suggest that flowcharts constructions
from [Ma3] can be interpreted as a constructive existence theorem, 
to produce what I call 
``an enriched programming method with unrestricted parallelism''
(cf. Definition 1.8.1). This seems to be very much in the spirit of [BaSt].

\smallskip

Finally, Section 3, using some ideas from quantum computation,
provides a Hopf renormalization scheme for the Halting Problem.

\smallskip

{\bf 0.2. Computability as a mathematical structure and its
interaction with  other mathematical structures.} 
Most of the constructions considered in this paper refer to (un)feasible algorithms
with infinite domains/ranges. Devising their natural quantitative characteristics   
and regularization schemes for them, one should keep in mind that
they can be roughly subdivided into two large blocks.

\smallskip

{\it BLOCK A.} This block consists of the inherent problems referring 
to an infinite constructive world $X$,
which depend {\it only on the class of
``admissible'' recursively  equivalent numberings of $X$,}
and which are the same for all infinite $X$. From the computational viewpoint,
any such $X$ can be identified with $\N$ (natural numbers) or $\Z$ (nonzero natural numbers).
\smallskip
A typical example of such an $X$ is some set of finite Bourbaki structures,
such words in a finite alphabet, or finite groups, or graphs,  or their descriptions, etc. 
In this context, one uses the Bourbaki description primarily in order to
define the class of admissible bijections (numberings) $\Z \to X$  in question: they must be
informally computable together with their inversions.
One aspect of Church's thesis consists in the statement 
that recursive functions will provide
an adequate notion of algorithmic processing of elements of $X$ 
whenever we can imagine an informal algorithm
producing the numbering.

\smallskip

Once it is decided that the role of the respective Bourbaki structure is over
as soon as the class of numberings is determined, one can
make explicit various secondary structures on $X$ 
that can be defined exclusively in terms of admissible numberings.
\smallskip
One  such structure is the  algebra of 
enumerable subsets of $X$: definition domains $D(f)$
of partial recursive functions. This family is stable wrt finite intersections and enumerable unions.
If one consider these sets modulo finite ones, one can prove interesting
results about simple sets, maximal sets etc.
For example, maximal $D(f)$ display a striking similarity 
to the holomorphy domains $V$
in the theory of complex analytic functions of $\ge 2$ variables: in both cases,
there are functions defined on $D(f)$, resp. $V$, that cannot be extended to
a larger domain.

\smallskip

The proviso ``modulo finite subsets'' can be very naturally formalized
by changing the Constructive Universe $\c$ described in Sec. 1: 
simply consider the largest quotient of $\c$
making invertible those morphisms (computable maps) $X\to Y$ that become
computably invertible after the restriction to some subsets of $X$ and $Y$ with finite complements.
\smallskip

Another such structure is the class of Kolmogorov's orderings:
total orders on $X$ defined by increasing Kolmogorov complexity with respect to
various optimal enumerations.
\smallskip

Such orderings  are not computable,
but with respect to them all recursive functions,
including admissible numberings,  become functions of linearly bounded growth.
I discuss this feature from the renormalization viewpoint in Sec. 3.

\smallskip

{\it BLOCK B.} This group consists of problems about 
interaction of computability with other Bourbaki sructures
on $X$. An elementary example is the
 embedding $X=\bold{Q}\subset \bold{R}$ used in the theory
of computable rational approximations to real numbers.
\smallskip
In this block, the Diophantine representability of enumerable subsets of $\N$
was the greatest discovery (Davis, Robinson, Putnam, Matiyasevich).

\smallskip

A very interesting and unexpected example of such interaction
was elaborated in the work of A.~Nabutovsky and S.~Weinberger,
(cf. [NaWe]),   who have shown that the computational complexity
can be used to display a highly irregular landscape
of minima of natural differential--geometric functionals
on the space of Riemannian metrics modulo diffeomorphisms.  

\smallskip

Since path integration over such a space is one of the key tools of
quantum gravity, this can become an important next meeting space
between  renormalization and computation. 

\smallskip

{\it Acknowledgements.} I am very grateful to Cristian Calude, Leonid Levin,
Mike Stay, Noson Yanofsky, who sent their remarks and suggestions
incorporated in this draft of the article.

\bigskip

\centerline{\bf 1. Enriched programming methods with typing and parallelism:}
\smallskip
\centerline{\bf a categorical approach}

\bigskip

{\bf 1.1. Preliminary remarks.}   We denote by $\N$ the set of natural numbers, and 
by $\Z$ that of positive natural numbers. Usually $\N$ is taken as the basic
set on which recursive functions are defined; in [Ma1] I used $\Z$ having found
it more convenient in a Diophantine context.

\smallskip

As a Bourbaki structure, both $\Z$ and $\N$ are here (isomorphic) totally ordered sets,
with a minimal element 1 (resp. 0) and successor function $\r{suc} (x)$:
``the smallest $y$ such that $y>x$'', or  $\r{suc} (x)=x+1$ in the standard notation.
In a sense, this is the minimal structure needed to define the set of partial recursive functions
that are partial maps  $\Z\to \Z$ or, more generally $(\Z )^a\to (\Z )^c$. 
The remaining components of the definition (see e.~g.~[Ma1], V.2.1--2.4)
are just the standard category--theoretic constructions in a fixed monoidal
category of sets $(ParSets, \times )$ with partial maps as morphisms and  cartesian product:
cf.~ [Ma3], 3.7. 
\smallskip

However this total order structure {\it is not} invariant with respect to the 
structure that we will define below, in the sense that it is not
preserved under the automorphisms of this structure. For this reason, we avoid
one of the standard categorifications of computation theory in which
$\N$ is replaced by a ``natural numbers object'' $\Cal{N}$ (of an abstract
category), endowed with
a morphism $\r{suc}_{\Cal{N}}:\,\Cal{N}\to\Cal{N}$: this categorification 
unduly stresses the role of this total order and iteration 
$\r{suc}_{\Cal{N}}\circ\dots\circ\r{suc}_{\Cal{N}}$  related to it.

\smallskip

Instead, we adopt the version advocated in [Ma2], that of a subcategory $\c$ of $ParSets$ 
 called {\it Constructive Universe}. In 1.2--1.8 below
 I collect the relevant formal definitions. Informal comments
 are relegated to 1.10.

\medskip

{\bf 1.2. Objects.} Objects of $\c$ will be called
constructive worlds.

\smallskip

{\bf 1.2.1. Definition.} {\it  A constructive world $X$ is either a finite set,
or an infinite set endowed with a nonempty  set $Num(X)$ of bijections $\nu :\Z\to X,$
called admissible numberings, satisfying  the following conditions:

\smallskip

(i) If $\nu_1,\nu_2\in Num(X)$, then $\nu_2^{-1}\circ\nu_1$ is a total recursive bijection.

\smallskip

(ii) If $\nu\in Num(X)$ and  $f:\,\Z\to \Z$ is a total recursive bijection,
then $\nu\circ f\in Num(X)$.

\smallskip

Elements of the constructive world $X$ are called
constructive objects of the type $X$.
}

\medskip

\medskip

{\bf 1.3. Morphisms.} Let $X,Y$ be two constructive worlds.
Morphisms $X\to Y$ are induced by partial recursive maps
on their structure numberings. More precisely:

\smallskip

{\bf 1.3.1. Definition.} {\it A morphism $X\to Y$ is partial map
$f:\,D(f)\to Y$, where $D(f)\subset X$ a subset (possibly empty)
satisfying the following conditions:

\smallskip

(i) If $X$ be infinite, $Y$ is finite, then for one (equivalently, any)
admissible numbering $\nu :\,\Z\to X$ and any $y\in Y$,
the set $\nu^{-1}(f^{-1}(y))$ is recursive enumerable.

\smallskip

(ii) If $X$ and $Y$ are infinite, then for one pair (equivalently, any pair)
admissible numberings $\nu_X:\,\Z\to X$, $\nu_Y:\,\Z\to Y$,
the partial map $\nu_Y^{-1}\circ f\circ \nu_X:\,\Z\to\Z$
is a partial recursive function.

\smallskip

(iii) If $X$ is finite, $Y$ is infinite, any partial map is a morphism.}
\smallskip

With the standard composition of partial maps, constructive worlds
form a category $\c$, which we will call {\it Constructive Universe.}
The set of morphisms $X\to Y$ will be denoted $\c (X,Y)$.
Its subcategory consisting of infinite constructive worlds
is equivalent to a very simple category consisting of one
object, say $\Z$, and partial recursive maps as morphisms.
In [He], such categories are called {\it isotypical} ones.

\smallskip

However, it is important to consider $\c$ as (bi)monoidal category,
with two symmetric monoidal structures $\times$ (direct product)
and $\coprod$ (coproduct, or disjoint union, see [He]), connected by the standard
coherence diagrams. 

\smallskip

These monoidal structures are induced by those in a small category of
(unstructured) sets in which our constructive worlds lie,
so that it suffices to specify some privileged numberings
of disjoint sums and direct products. Moreover, it suffices
to consider numberings that are bijective maps  
$\Z\to\coprod_{i=1}^m \Z$ and $\Z\to (\Z)^m$.

\smallskip

For   $\coprod_{i=1}^m \Z$, we simply assign to $m(k-1)+i$ 
the number $k$ in the $i$--th summand.

\smallskip

The cartesian product is more interesting, because there are
several numberings that become privileged in the context
of Kolmogorov complexity. {\it TO BE CONTINUED ...} 

\smallskip

We will generally assume that $\c$ is closed with respect to
the monoidal structures $\times$ and $\coprod$.

\medskip

{\bf 1.4. Constructive descriptions of morphisms.} Fix two
constructive worlds $X,Y$.
Since  the set of recursive maps $\c (X,Y)$ is not
a constructive world, we may try to replace it by descriptions.

\medskip

{\bf 1.4.1. Definition.}  {\it   A constructive world of descriptions  is a pair
$(P(X,Y), F)$, where $P(X,Y)$ is an object of $\c$, and $F:\, P(X,Y)\times X\to Y$
is a morphism in $\c$, satisfying the following condition.

\smallskip

Let $p\in P(X,Y)$. Denote by $f_p$ the partial map
$x\mapsto F(p,x)\in Y, x\in X$. Then each $f_p$ is a morphism $X\to Y$ in $\c$.
}

\smallskip
In other words, descriptions produce  a set theoretic map $P(X,Y)\to \c(X,Y)$
constructively depending on $(p,x)$.

\medskip

{\bf 1.4.2. Translations.} Let $(P(X,Y),F)$ and $(Q(X,Y),G)$ be two
constructive worlds of descriptions. A {\it translation} (or {\it compilation}) 
{\it method}
$$
trans_{P,Q}:\, P(X,Y)\to Q(X,Y)
$$
is an everywhere defined morphism in $\c$  such that for all $p\in P(X,Y)$,
$trans (p)\in Q(X,Y)$ defines the same morphism $f_p:\,X\to Y$.
In other words,
$$
G\circ (trans_{P,Q}\times id_X)= F.
$$

\smallskip

{\bf 1.4.3. Universal descriptions.} The world $(U(X,Y),W)$ is the world
of {\it universal descriptions,} if for any other world of descriptions
$(P(X,Y),F)$, there exists a translation morphism
$$
trans_{P,U}:\, P(X,Y)\to U(X,Y)
$$ 

In particular, it can compute any (semi)computable function,
in the sense that the family of maps $f_u:\,X\to Y,\, u\in U(X,Y)$
contain all morphisms in $\c$.

\medskip
{\bf 1.4.4. Complements.} Among various constructive 
worlds of descriptions $P(X,Y)$ 
there exist ones with better properties than the general
definition allows to guess. The terminology for them
is rather unstable. We will say sometimes (see 2.6 below) that
$P(X,Y)$ is {\it the base of a family} $f_p:\,X\to Y$, $p\in P(X,Y)$
of partial functions.

\smallskip

We will review below
some relevant definitions and existence theorems.
We consider separately four cases.

\smallskip
(i) If $X,Y$ are infinite, we may assume 
without losing generality that $X=Y=\Z$.
\smallskip

H.~Rogers in [Ro]  calls such a world of descriptions $U$  {\it semi--effective},
if it computes all morphisms (partial recursive functions), and he calls $U$ {\it fully effective},
if it is universal in the sense of 1.4.5.

\smallskip
 An easy construction in [Ro] (example following Definition 3) shows that 
 there are semi--effective descriptions that are not fully effective. 
 A universal description world  $U(X,Y)=\Z$
 (or rather $\N$) is also called a {\it G\"odel numbering}
in [Ro].

\smallskip

The main Theorem of [Ro] in our language implies that {\it for any
two universal description worlds (for infinite $X$, $Y$),
there exist two mutually inverse translation isomorphisms between them:
total recursive bijections, compatible with functions that these 
descriptions compute.}

\smallskip

C.~P.~Schnorr in [Sch] considerably strengthens this result.
Namely, he calls a universal description world $(U=\Z,W)$  {\it an optimal
G\"odel numbering}, if for any other world of descriptions $(P=\Z, F)$, 
there exists a translation morphism $t:\,P\to U$ which is a
linearly bounded function $\Z\to\Z$. We will call such a description world simply {\it optimal} one.

\smallskip

Schnorr then proves that  {\it optimal descriptions exist, and for any
two optimal description worlds (with infinite $X$, $Y$),
there exist two mutually inverse linearly bounded translation isomorphisms between them}.

\smallskip

Similar results hold in the case when only one of the worlds $X,Y$ is infinite.
The general situation can be reduced to the case when the relevant finite
world is one--element set; for simplicity, we will consider only this case. 

\smallskip

(ii) The case $X=\{*\}$ is truly exceptional in the following sense:
$\c (\{*\},Y)$ , that is ``$Y$--valued recursive functions of zero variables'',
can be canonically identified with the set $Y$ 
and thus {\it it is} a constructive world. Nevertheless,
the notions and main results of Rogers and Schnorr
are applicable to this case as well and lead to a strengthening
of the notion of Kolmogorov optimal enumeration of constructive objects
of a given type.

\smallskip

(iii) In the case $Y=\{*\}$,   $\c (X,\{*\})$ can be naturally identified
with the set of all enumerable subsets of $X$, domains of
partial recursive functions with one value. G\"odel and optimal
numberings of enumerable subsets also can be easily defined, and again
the Rogers and Schnorr theorems are valid for them.

\smallskip
(iv) Finally,  when both $X,Y$ are finite, the useful structurizations of
descriptions are those  of Boolean polynomials, circuits, etc.
Many complexity problems are centered around polynomial time computations.
Cf. [Ma2] for an introduction, that is close in style to this paper.

\medskip

{\bf 1.5. Enrichments of $\c$ over itself and programming methods.}
There is a well--known general notion of a category $C$ enriched over
a monoidal category $(M,\otimes , I)$ where $I$ is an identity object.

\smallskip

Below, we will consider enrichments of $\c$  over $(\c,\times ,I)$ where
$I$ is a fixed one--element constructive world. Products of empty families,  
such as $\N^a$ for $a=0$, are interpreted as $I$.

\smallskip

According to the general pattern, such an enrichment must consist of the
following data.

\smallskip

a) For each pair of constructive worlds $X,Y$, an "object of morphisms" $P(X,Y)\in \r{Ob} \c$.

\smallskip

b) For each triple of constructive worlds $X,Y,Z$, a ``composition morphism''
in $\c$:
$$
\circ :\, P (Y,Z)\times P (X,Y)\to P (X,Z)
\eqno(1.1)
$$

c) For each object $X$ of $\c$, an identity morphism
$$
\r{id}_X:\,I\to P (X,X).
\eqno(1.2)
$$

\smallskip

The standard axioms for morphisms in a category translate into
requirements of commutativity of three classes of  diagrams 
in $\c$ expressing properties of associativity of enriched
composition, and left and right identities.

\medskip

{\bf 1.5.1. Definition.} {\it An enrichment of $\c$  over $(\c ,\times , I)$ as above
is called an enriched programming method and denoted $\c_P$, if the following additional data are given and axioms satisfied:

\smallskip

For each pair of constructive worlds $X,Y$, a morphism in $\c$ is given
$$
F_{X,Y}:\, P(X,Y)\times X\to Y
$$
such that $P(X,Y)$ becomes a constructive world of descriptions in the sense of
Definition 1.4.1.  Thus, we have a family of set--theoretic maps
$$
\Phi_{X,Y}: \, P(X,Y)\to \c (X,Y): \quad p\mapsto f_p.
\eqno(1.3)
$$
We  will say that $p$ is a description, or a program, computing $f_p$.

\smallskip

Moreover, the following axioms must be satisfied:

\smallskip

(i)  Morphisms (1.1) and (1.2) must be everywhere defined (total recursive) maps.
(For (1.2), this means simply that they are non-empty maps).

\smallskip

(ii) Compositions (1.1) must be compatible with the compositions
of  morphisms in $\c$:
$$
f_{p\circ q} =f_p\circ f_q.
$$

(iii) Element $\r{id}_X(I)\in P(X,X)$ must be a description of the 
``copying'' program: mapping $x$ to $x$.

\smallskip

(The latter should not be mixed with  the ``cloning'' program computing the dia\-gonal map $X\to X\times X, x\mapsto (x,x)$.)}
\medskip

Informally speaking, we have a functor
$$
\Phi_P :\, \c_P \to \c
\eqno(1.4)
$$
identical on objects and mapping a program to the function that this program computes.

\smallskip

The most important feature of this formalism consists in 
an explicit and systematic inclusion of composition of programs
into our formalism: this is a key requirement for all Hopf
algebra renormalization schemes.

\smallskip

Finally, a remark on terminology: we use the word {\it program}
as a synonym of  {\it ``description of a method to compute a given
function''}. Input of such a program is a specific value of the argument
of this function, output is the value of the function.
This praxis should not be confused with the one 
used in the theory of Turing machines, where
programs are often understood as our inputs: the initial
binary string on the tape.

\medskip

{\bf 1.6. A (uni)versal enrichment.}  An enrichment $\c_U$ as above
is called {\it (uni)versal} one, if  programs from $U$
compute {\it all} partial recursive maps, and moreover, if  for each $\c_P$, there is a functor
between enriched categories 
$$
\Psi :\, \c_P \to \c_U
\eqno(1.5)
$$
identical on objects, with total recursive maps
$$
\Psi_{X,Y)}:\, P(X,Y)\to U(X,Y)
\eqno(1.6)
$$
such that $\Psi_{X,Y}(p)$ for each $p$
computes the same function as $p$. In other words, we have 
$$
\Phi_U\circ\Psi =\Phi_P.
$$
Intuitively, we want the following properties of $U$ as a programming method,
that are somewhat stronger than those in 1.4.

\medskip

a) It can compute any (semi)computable function.

\smallskip

b) For each other programming method $P$, there must exist
a computable (on the world of $P$--programs)  {\it translation}  of $P$--programs 
into $U$--programs,
computing the same functions.

\smallskip

c)  The translation must be compatible with composition of programs
and copying/identity programs (functorality of $\Psi$).

\medskip

{\bf 1.7. Coproducts and typing.}  Coproducts (disjoint sums) admit the most
straightforward interpretation in the contexts, where computer scientists
speak about {\it typing}. In the simplest situation,
a program $p\in P(X\coprod Y, Z)$ accepts inputs of either type $X$
or  $Y$, and produces outputs of type $Z$. 

\smallskip

Iterated application of similar interpretations, as far as I can judge,
can be used in all contexts where the notion of typing is essential.

\medskip

{\bf 1.8. Products and parallelism.} Let $X_i,Y_i, i=1,\dots, n$
be  constructive worlds. Then the map
$$
\pi :\,\c(X_1, Y_1)\times \c (X_2,Y_2)\times\dots \times  \c(X_n,Y_n)\to \c(X_1\times X_2
\times \dots \times X_n, Y_1\times Y_2\times \dots \times Y_n) ,
$$
$$
\pi (f_1,\dots , f_n) (x_1,\dots ,x_n) := (f_1(x_1),\dots ,f_n(x_n))
\eqno(1.7)
$$
defines families of computations with independent inputs/outputs
that can be implemented parallely. This can be generalized to
programming methods as follows.

\smallskip

{\bf 1.8.1. Definition.} {\it Let  $P$ be an enriched programming method.
We will say that $P$ admits unrestricted parallelism, if the following additional
structure is given.

\smallskip

Let $\{X_i\}$, $\{Y_i\}$, $i=1,\dots ,n$, be any finite family $\sigma$ of constructive worlds.
We must be given maps
$$
\pi_{\sigma} :\,  \prod_{i=1}^n P(X_i, Y_i) \to P(\prod_{i=1}^nX_i,\prod_{i=1}^nY_i) ,
\eqno(1.8)
$$
that are lifts of (1.7). These maps must be equivariant with respect to the natural 
action of the symmetric group $S_n$ permuting subscripts at both sides. }

\medskip

{\bf 1.9. Basic example: flowcharts.} Flowcharts defined in [Ma3] form
a convenient context for constructing enriched progamming methods with unrestricted parallelism. One example of  of such a method is
the world $P$, computing primitive recursive functions, 
described im [Ma3], Definition 2.11. This definition uses
the ideas of N.~Yanofsky paper [Ya].

\smallskip

Additional work remains to  be done in order to produce a manageable
construction of an universal enrichment with unrestricted parallelism.

\medskip

{\bf 1.10. Comments: constructive worlds and
admissible numberings.}  Technically speaking,
conditions  (i) and (ii) of Definition 1.2.1 together mean that 
$Num(X )$ forms a principal 
homogeneous space over the group of total recursive permutations of $\Z$ which we 
may denote $S_{\Z,rec}$ , as an infinite analog of Þnite symmetric groups. 
Hence the whole $Num(X )$ can be reconstructed from any one 
numbering in this set. Usually there 
are some ÒsimplestÓ, or ÒprivilegedÓ numberings, with which one mostly works, 
such as the numbering of binary words $w  \in \{0, 1\}^{\N}$ 
used in [CalSt1]: ${bin}^{-1} : w\mapsto \overline{1w}$, 
where the line over a binary word means that it should be treated as a natural 
number given by its binary digits.

\smallskip 

The idea of privileged numberings is essential especially when one deals with 
polynomial time, or more general ``feasible'' computations. 
In order to accommodate this idea, we can strengthen Definition 1.2.1 in the 
following way. Consider the smaller ``symmetric group'' $S_{\Z,pol}$ of total recursive 
permutations that are polynomial time computable together with their inverses. 
Define the structure of {\it polynomial time constructive world} $X$ 
by a set $Numpol (X )$ 
forming a principal homogeneous space over the group $S_{\Z,pol}$ . 

\smallskip

An (often implicit) part of contemporary philosophy
around  Turing's Thesis consists in postulating
that whenever we can informally speak about algorithms and (semi) computable
maps between two constructive worlds, we always can 
produce in the context of this discussion admissible numberings that
are informally algorithmic and transform informal
(semi)computable maps into (partial) recursive functions. 

\smallskip

In any case, starting with such a class of numberings of $X$, we want to stress that we 
study notions that are either invariant, or behave in a controlled way under the 
action of the  group  $S_{X,rec}$. 
\smallskip
We add a few more remarks. 

\smallskip
Sometimes, a constructive world $X$ itself is an unstructured set, in the sense that 
the only relevant structure on it is given by its set of admissible numberings. The 
typical example is a world $A$ that in further constructions may serve as an alphabet. 
In this case, a numbering deÞning the whole $Num(A)$ usually is introduced ad hoc. 
But in most applications, $X$ itself consists of certain sets (often finite and/or 
considered only up to isomorphism, or even organized into a category) endowed 
with a certain fixed Bourbaki structure, such as: 
\smallskip
(a) Finite words in an alphabet $A.$ 

\smallskip
(b) Finite graphs up to an isomorphism. 

\smallskip
(c) Finite groups. 
\smallskip
In such cases, the privileged numberings (``encodings'')
generating the whole $Num (X)$
 are supposed to interact with this structure in such 
a way that the number of a constructive object can be ``algorithmically calculated''
when we know this object as an instance of this structure, and vice versa, this 
instance must be algorithmically reconstructible from the number.
It is only very rarely that such encodings can translate well the
basic relations, composition laws etc., involved in the definition
of the structure. This is one reason to formulate models of computation 
directly in terms of this structure: one can recall Church's
lambda--calculus, Kolmogorov--Uspenski's graphs and G\'acs--Levin 
causal nets (cf. [GaLe]). 

\smallskip

On the other hand, such a simple task as the choice of a privileged
numbering of $\Z\times \Z$ (upon which one of the monoidal structures
on $\c$ is based) can lead to quite interesting constructions when this
choice is related, for example, with complexity
estimates: see our discussion in 2.7--2.10 below using
L.~Levin's norms  for definitions of such numberings.

\smallskip

Notice that instances of the constructive worlds of causal nets 
studied in [GaLe] are themselves {\it categories},
and the study of interaction of computability
with symmetries of the respective constructive objects
reveals interesting new phenomena. 

\smallskip

We want to argue that our categorical framework suggests other 
possibilities to avoid 
too close attention to the elementary steps of computation. 
In particular, the ``categorical Church Thesis" admits a wonderfully
succinct expression: 

\smallskip
{\it The category $\c$ is defined uniquely up to equivalence}.

\smallskip
    
An important complication and variation of the theme
of admissible numberings arises, when a structure 
$S$ that we 
want to treat ``constructively'' is thus imposed on eventually inÞnite sets.
In such cases  the 
relevant constructive objects are often not the structures themselves,
but their finite  {\it descriptions}: a group might be given by generators and
relations, an affine scheme over $\bold{Z}$ by its equations etc.

\smallskip

The usual complication with descriptions is that many descriptions can produce 
one and the same (or canonically isomorphic) Bourbaki structure, and the relevant 
equivalence relation on the set of descriptions can be {\it undecidable,} or even 
{\it not recursively enumerable.} This is precisely the case of the structure 
constituted by recursive functions themselves, which is our main motivation 
for introducing enrichments as in Definitions 1.5.1 and 1.8.1. 

\bigskip

\centerline{\bf 2. Cut--off regularization and Anytime Algorithms}

\medskip

{\bf 2.1. Cut--off regularization.} In Quantum Field Theory, cut--off regularization
schemes have the following typical structure. The relevant Feynman integrals, say, in
momentum space, may diverge when momentum becomes large (resp. small).
In this case, the formal integral in question $I$ is replaced by the finite integral $I_P$
taken over momenta $p\le P:=p_{cutoff}$ (resp. $p\ge P:=p_{cutoff}$).
The behavior of $I_P$ as $P\to \infty$ (resp. $P\to 0$) is then studied, and 
physical information is extracted from the behavior of the polar part,
or regular part, of  $I_P$. 

\smallskip

In computer science based upon Turing machines and/or
recursive functions, the natural  ``divergence'' occurs in space--time:
a computation uses discrete memory (space) and runtime.
A typical example of such a divergence is the infinite 
runtime of a Turing machine computing a partial
recursive function $f$ at an input (program) $x$ which is outside
the definition domain of $f$.
\smallskip

Application--oriented computer scientists, of course,
 recognize the practical necessity
of time cut--offs, accompanied by sober estimates of
quality of outputs. Systematic work on this problem 
resulted in the notion of "Anytime Algorithms",
cf. [GrZi]. The usefulness of composition and exploiting parallelism
was stressed in [RuZi].

\smallskip

In a stimulating paper [CalSt1], Ch.~Calude and M.~Stay addressed the problem
of  cut--off  of runtime theoretically, and designed  meaningful 
quantitative characteristics of such a cut--off. 

\smallskip

More precisely, let $f$ be a partial recursive function (``a morphism of
constructive worlds $X\to Y$'') as above,
and $F$ its description, a program calculating it.  {\it Computation time (or runtime)} is
another partial recursive function, with the same
domain $D(t(F))=D(f)\subset X$ and target $\Z$, whose precise definition depends 
on the details of the implied choice of our programming method.

\smallskip
For a Turing machine $F$, $t(F)(x)$  is the number of steps required
to halt and print $f(x)$ on the tape, if $x\in D(f)$. One can similarly
define another partial recursive function, 
``memory volume'' $m(F):\,X\to\Z$ such that
for $x\in D(f)$, $m(F)(x)$ is the minimal length of tape required 
to compute $f(x)$. Here $W$ is the constructive world of binary
words $\{0,1\}^{\Z}.$ Yet another partial recursive function, $s(F)$, is
the sum total of lengths of filled parts of the tape over all steps
of computation. Notice that the {\it settling function} of Soare ([So], Definition 8.2),
which is essentially $\r{max}\,\{ t(F)(y)\,|\,y<x, y\in D(f)\}$, generally is not
partial recursive, but some of the inequalities stated below,
such as (2.3), are valid for it as well.
\smallskip
One can define natural  analogs of functions $t(F)$, $m(F)$, and $s(F)$
for rather general normal programming methods $F$, discussed in  [Ma2]
and Chapter IX of the new edition of [Ma1]. 

\smallskip

Returning to [CalSt1], we will first of all show that some of the basic results of
that paper related to cut--offs, admit a straightforward reformulation 
in such a way that they become  applicable and true for
{\it any} partial recursive function, including, of course,
$t(F)$, $m(F)$, and $s(F)$.

\smallskip

This naturally raises a question, what is so specific about
$t(F)$, $m(F)$, and $s(F)$. We will treat this question in 2.6 below in
the context of categorification developed in Sec.~1, and will show 
that this provides some meaningful insights about these
measures of processes of computation. 

\medskip

{\bf 2.2. Complexity.} I will first recall the definition and properties
of the Kolmogorov (``exponential'', or ``program'') complexity $C_u:\,\Z \to \Z$,
cf. [Ma1], VI.9. In [CalSt1] it is called {\it the natural complexity} and denoted
$\nabla_U$ or simply $\nabla$. 

\smallskip

This complexity measure is defined with respect to a
partial recursive function $u:\,\Z\to \Z$ which is surjective:
$$
C_u(x):=\r{min}\,\{y\,|\,u(y)=x\}.
$$
This function $u$ is  an arbitrary element of the set
of
{\it Kolmogorov, or G\"odel,  optimal} functions, representatives of which can be
effectively constructed: cf. [Ro], [Sch] and [Ma1].  Optimality
implies that for any other partial recursive
$v:\,\Z\to\Z$, there exists a constant $c_{u,v}>0$ such that
for all $x$, $C_u(x)\le c_{u,v} C_v(x)$. (The right hand side is interpreted
as $\infty$, if $x$ is not in the range of $v$). 

\smallskip

It follows that another choice of optimal function 
replaces $C_u$ by a function $C_{u^\prime} = 2^{O(1)}C_u$.
We will say that two such functions belong to the same
{\it bounded equivalence} class.

\smallskip

Moreover, as we have discussed above, the
definition of the complexity of integers can be extended
to the definition of complexity of  partial recursive functions
of any fixed number of variables $m$, as in VI.9.1 of [Ma1].
This requires a choice of Kolmogorov optimal 
recursive function of $m+1$ variables. 
We have then the following simple result (omitting
the subscripts  at $C$ specifying the choices of optimal families,
and denoting by $c$ with subscripts various constants depending
on these choices as well):

\medskip

{\bf 2.2.1. Proposition.} {\it For any partial recursive function $f:\, \Z\to \Z$
and $x\in D(f)$
we have
$$
C(f(x))\le c_f C(x) \le c^{\prime}_f x,
\eqno(2.1)
$$
If $f$ and $x\in D(f)$ are allowed to vary, we have }
$$
C(f(x))\le c\, C(f) C(x)\,\r{log}\,(C(f)C(x))\,.
\eqno(2.2)
$$

\smallskip

In particular, if $f$ is a total recursive permutation, then complexities
of $x$ and $f(x)$ are bounded equivalent. It follows that
we can define the complexity function, up to bounded equivalence,  $C:\,X\to \Z$ for any
infinite constructive world $X$: choose an admissible
numbering  $\nu :\,\Z\to X$ and put $C(x):=C_u(\nu^{-1}(x))$
for some optimal $u$.

\medskip

{\bf 2.3. Runtimes according to [CalSt1].} Proposition 2.2.1 is 
a special case of Proposition 9.6 in [Ma1], VI.9.
In turn, it implies as special cases the inequality (2) and Theorem 4
of [CalSt1].

\smallskip

 In order to see this, one has  simply to compare terminology and
notation.

\smallskip

[CalSt1] deals with the complexity $C$ (their $\nabla =\nabla_U$)
of binary words, that reduces to the complexity of integers
via the admissible numbering denoted ${bin}$ in [CalSt1].
It is defined via a "universal Turing machine" $U$,
which in our language is a programming method computing one of 
the Kolmogorov optimal functions $u$. Consider the partial recursive function
$x\mapsto t(U)(x)$: runtime of $U$ at the argument $x\in \Z$.  
The inequality (2) of [CalSt1] in our notation
can be rewritten as
$$
C(t(U)(x)) \le cx
\eqno(2.3)
$$
which is our (2.1) for $f=t(U)$. The same inequality is valid for $m(U)$,
$s(U)$, but also for $t(F)$, $m(F)$, $s(F)$ for any $F$,
and even for Soare's settling functions: see section 2.1 above.

\medskip

{\bf 2.4.1. Growth of recursive functions and algorithmic randomness.}
The central argument of [CalSt1] is based upon two statements:

\smallskip

{\it a) The runtime of the Kolmogorov optimal program at  a point $x$
of its definition domain is either $\le cx^2$, or is not
``algorithmically random'' (Theorem 5 of [CalSt1]).

\smallskip

b) ``Algorithmically random'' integers  have density zero for a class
of computable probability distributions.}

\smallskip

This last statement  justifies the time cut--off prescription which is the main result
of [CalSt1]:
\smallskip
{\it if the computation on the input $x$  did not halt after $cx^2$ Turing steps,
stop it, decide that the function is not determined at $x$, and proceed to $x+1$.}

\smallskip

Proposition 2.5.1 below somewhat generalizes the statement a).

\medskip

{\bf 2.5. Randomness and growth.} Consider a pair of functions
$\varphi ,\psi :\, \R_{>0} \to \R_{>0}$ satisfying the following conditions:

\smallskip

a) $\varphi (x)$ and $\dfrac{x}{\varphi (x)}$ are strictly increasing
starting with a certain $x_0$  and tend to infinity as $x\to\infty$.

\smallskip

b) $\psi (x)$ and $\dfrac{\psi (x)}{x\varphi (\psi (x))}$ are increasing and tend 
to infinity as $x\to\infty$.

\smallskip

The simplest examples are $\varphi (x)=\r{log} (x+2)$, $\psi (x)=(x+1)^{1+\varepsilon}$,
$\varepsilon >0$.

\smallskip
In our context, $\varphi$  will play the role of  a  ``randomness scale''.
Call $x\in \Z$ {\it algorithmically $\varphi$--random,} if $C(x)>x/\varphi (x).$
The second function $\psi$ will then play the role of associated growth scale.
\medskip

{\bf 2.5.1. Proposition.} {\it  Let $f$ be a partial recursive function.
Then for all sufficiently large $x$ exactly one of the following alternatives holds:

\smallskip

(i) $x\in D(f)$, and $f(x)\le \psi (x)$.

\smallskip

(ii) $x\notin D(f)$.

\smallskip

(iii) $x\in D(f)$, and $f(x)$ is not algorithmically $\varphi$--random.}

\medskip

{\bf Proof.} We must only check that if $x\in D(f)$ and $f(x)>\psi (x)$,
then $f(x)$ is not algorithmically $\varphi$--random, that is
$$
C(f(x))\le\frac{f(x)}{\varphi (f(x))}.
\eqno(2.4)
$$
In fact, in view of (2.1), 
$$
C(f(x))\le cx
\eqno(2.5)
$$ for some constant $c$ (depending on $u$ and $f$).
Furthermore, for sufficiently large $x$, in view of 2.5 b), we have
$$
cx\le \dfrac{\psi (x)}{\varphi (\psi (x))}\le \frac{f(x)}{\varphi (f(x))}.
\eqno (2.6)
$$    
Clearly, (2.5) and (2.6) imply (2.4).

\medskip

{\bf 2.6. Cost estimate functions.}  Since, as we argued, the
randomness/growth alternative holds for arbitrary recursive functions,
not only for runtimes and alike, we will briefly discuss specific properties
of runtimes, considered from the perspective of categorification,
explained in sec. 1.

\smallskip

Let $P$ be an enriched programming method, as in Definition 1.5.1.
We will say that a partial function $ \delta :\,P(X,Y)\times X\to \Z$ 
is {\it a cost estimate function}, if the following conditions are satisfied:
\medskip

{\it (i) $\delta$ is partial recursive (morphism in $\c$), and $D(\delta )=
\{(p, x)\,|\,x\in D(f_p)\}$.

\medskip

(ii) $\delta (p\circ q,x)=\delta (q,x) + \delta (p, q(x))$ whenever both sides
are defined.}

\medskip

The requirement (i) is natural, because the ``run--cost'' of computation
(time, maximum storage size) must be computable in terms
of cost increments required at each step. The requirement
(ii) then expresses the additivity of such increments.
We may, or may not, ascribe a non--zero cost to the
program calculating identical function
(``data transfer'').

\smallskip

Requirement (i), complemented by the requirement of the {\it decidability
of the graph of $\delta$}, constitute two axioms due to M.~Blum.
This latter property has  a clear intuitive meaning as well.

\smallskip

Finally, if our cost estimate function refers to time only, and we allow
the unrestricted parallelism, the following
property is natural. Using notation (1.8), we must have

\medskip

{\it (iii) $\delta (\pi_{\sigma}(p_1,\dots ,p_n),(x_1,\dots ,x_n))=
\r{max}\, (\delta (p_1,x_1),\dots ,\delta (p_n,x_n)).$}

\medskip

{\bf 2.7. Constants related to Kolmogorov complexity estimates.} Since inequalities (2.1), (2.2), and their extensions are often useful, we will say a few words
about their computability.

\smallskip

As in [Ma3], VI.9, we call partial maps $(\Z)^m\to \Z$, $m\ge 0$,
{\it m--functions.} A Kolmogorov optimal family of $m$--functions
$u(x_1,\dots ,x_m;k)$, $k\in \Z$, is produced from two inputs:

\smallskip

(a) A  fully effective (in the sense of Rogers, cf. 1.4.4)
family of $(m+1)$--functions $U$.

\smallskip

(b) A recursive embedding $\theta :\, \Z\times\Z\to \Z$
with decidable image, satisfying a linear growth condition
$$
\theta (k,j)\le k\cdot \varphi (j)
\eqno(2.7)
$$
where $\varphi :\,\Z\to\Z$ an appropriate function.

\smallskip

Having made these choices, we put
$$
u(x_1,\dots ,x_m;k):= U(x_1,\dots ,x_m; \theta^{-1}(k)).
\eqno(2.8)
$$
Then for any other family  $v$ of $m$--functions $v$ with base $\Z$
and each $m$--function $f$,
we have the inequality
$$
C_u(f)\le c_{u,v}C_v(f),
\eqno(2.9)
$$
with
$$
c_{u,v}:= \varphi (C_U(v)).
\eqno (2.10)
$$
(cf. [Ma1], VI.9.4).

\smallskip
Clearly, (2.1) is a special case of (2.9). An effective estimate
of (2.10) from above
will be assured, if $\varphi$ is computable and increasing, and if,
knowing a $P$--description of $v$, we can find some member of the family $U$ coinciding with
$v$. The latter, in turn, is automatic, if $P$ is supplied with
a translation morphism $trans_{P,U}$. 

\smallskip

We will now discuss numberings $\theta$.

\medskip

{\bf 2.8. Slowly growing numberings.} Let $R=(R_k\,|\,k\in\Z)$ 
be a sequence of positive  numbers tending to infinity
with $k$. For $M\in \Z$, put
$$
V_R(M):=\{(k,l)\in (\Z)^2 \,|\, kR_l\le M\}.
\eqno(2.11)
$$
Clearly,
$$
\r{card}\,V_R(M)\le \sum_{l=1}^{\infty} \left[\frac{M}{R_l}\right] < \infty\,,
\eqno(2.12)
$$
where $[a]$ denotes the integral part of $a$.

\smallskip

We have 
$$
V_R(M)\subset V_R(M+1),\  (\Z)^2=\cup_{M=1}^{\infty} V_R(M).
$$

Therefore we can define a bijection $N_R:\,\Z\to (\Z)^2$ in the following
way: $N_R(k,l)$ will be the number of $(k,l)$ in the total ordering
$<_{R}$ of $(\Z)^2$ determined inductively by the following rule:
$(i,j)<_R(k,l)$ iff one of the following alternatives holds: 
\smallskip
(a) $iR_j<kR_l$;
\smallskip
(b) $iR_j=kR_l$ and $j<l$;
\medskip

{\bf 2.9. Proposition.} {\it The numbering $N_R$ is well defined
and has the following property: all elements of $V_R(M+1)\setminus 
V_R(M)$ have strictly larger numbers than those of $V_R(M)$. Moreover:
\smallskip

(i)  If each $R_l$ is rational, or computable from above, then
$N_R$ is computable (total recursive).

\smallskip

(ii) If the series $\sum_{l=1}^{\infty} R_l^{-1}$ converges and its sum is bounded by a constant $c$, then
$$
N_R(k,l)\le c(kR_l+1).
\eqno(2.13)
$$
\smallskip

(iii) If the series $\sum_{l=1}^{\infty} R_l^{-1}$ diverges, and  
$$
\sum_{l=1}^{M} R_l^{-1}\le F(M)
\eqno(2.14)
$$
for a certain increasing function $F=F_R$, then
$$
N_R(k,l)\le (kR_l+1) F(kR_l +1).
\eqno(2.155)
$$
}

\medskip

{\bf Proof.} The first statements are an easy exercise.
For (2.13) and (2.14), notice that if $M$ is the minimal value for which
$(k,l)\in V_R(M)$, we have  $M-1<kR_l\le M$ and
$$
N_R(k,l) \le \r{card}\,V_R(M),
$$
and in the case (ii) we have from (2.12) 
$$
\r{card}\,V_R(M)\le \sum_{m=1}^{\infty} MR_m^{-1} \le c(kR_l+1).
$$

Similarly, in the case (iii) we have
$$
\r{card}\,V_R(M)\le M\sum_{m=1}^M R_m^{-1}\le (kR_l+1)F(kR_l+1)
$$
\smallskip

{\bf 2.10. L.~Levin's norms.}  From (2.13) one sees that any sequence 
$\{R_l\}$ with
converging  $\sum_lR_l^{-1}$ can be used 
in order to construct the bijection
$\Z\times \Z\to \Z$, $(k,l)\mapsto N_R(k,l)$
linearly growing wrt $k$. Assume that it is computable
and therefore can play the role of $\theta$ in (2.7) (b).

\smallskip

In this case, for any integer $M$ the set $V_R(M)$ must be decidable.
It follows that for any $l$, the set of rational numbers
$k/M\le r_l:=R_l^{-1}$ is decidable. 

\smallskip

Even if we  weaken the last condition, requiring only recursivity
of the set $k/M\le r_l$ (i.~e.~ asking each $r_l$ to be computable from below), the convergence of $\sum_l r_l$
implies  that there is a universal upper bound (up to a constant) for such 
$r_l$. Namely, let $CP$ be  the {\it prefix Kolmogorov complexity}
on $\Z$ defined with the help of a certain optimal prefix
enumeration.

\medskip

{\bf 2.10.1. Proposition.} ([Le]). {\it  For any sequence of 
computable from below numbers $r_l$ with convergent
$\sum_l r_l$, there exists a constant $c$ such that for all $l$,
$r_l\le c\cdot CP(l)^{-1}$ }

\medskip

More generally, L.~Levin constructs in this way a hierarchy
of complexity measures associated with  a class of abstract
norms, functionals on sequences computable from below.

\bigskip

\centerline{\bf 3. Regularization and renormalization of the Halting Problem}

\medskip

{\bf 3.1. Introduction.} In this section, we devise simple  
regularization/renormali\-zation  schemes
tailored to fit the halting problem. The general structure
of such a scheme is sketched in [Ma3], subsection 0.2.
It involves the following components.

\smallskip

{\it (a) Deforming the Halting Problem.}  At this step, we
transform the problem of recognizing, whether
a number $k\in \Z$ belongs to the definition
domain $D(f)$ of a partial recursive function $f$,
to the problem, whether an analytic function $\Phi (k,f;z)$
of a complex parameter $z$  has a singularity (in our case, a pole) at $z=1$.

\smallskip

In fact, using an idea from quantum computing, we may  reduce
the case of arbitrary $f$ to the case of a partial recursive 
permutation $\sigma=\sigma_f:\, D(\sigma ) \to D(\sigma )$ of its  definition domain,
and construct $\Phi (k,\sigma ;z)$ for such permutations.
This reduction is described in [Ma3], subsections 3.6 -- 3.8.

\medskip

{\it (b)  Choosing a minimal subtraction algebra.} Our choice of
an appropriate minimal subtraction algebra (see the definition in [Ma3], 4.2) is based on the established properties
of functions $\Phi (k,\sigma ;z)$: cf. Proposition 3.5 below.

\smallskip

Namely, let $\Cal{A}_+$ be the algebra of analytic functions
in $|z|<1$, continuous at $|z|=1$. It is a unital algebra; we endow it
with augmentation $\varepsilon_{\Cal{A}}:\, \Phi (z)\mapsto \Phi (1).$ 
\smallskip
Put
$\Cal{A}_-:= (1-z)^{-1}\C [(1-z)^{-1}].$ Finally, let 
$\Cal{A}:=\Cal{A}_+  \oplus \Cal{A}_-$.

\smallskip

Now we can use Theorem 4.4.1 of [Ma3] in renormalization schemes,
involving a connected filtered  Hopf algebra $\Cal{H}$
(cf. [E-FMan], sec. 2.5, Theorem 1, and [Ma3],  4.1). It remains to indicate, which Hopf algebras
and their  $\Cal{A}$--characters will be involved  in this game.

\medskip

{\it (c) Hopf algebra of an enriched programming method.} A class
of such algebras is described in [Ma3], subsections 3.3 -- 3.4.
This construction explicitly refers to flowcharts, however,
it can be readily modified and generalized to enriched
programming methods $P$ in the sense of Definition 1.5.1.
\smallskip
Basically, $\H=\H_P$ is the symmetric algebra, spanned by
isomorphism classes $[p]$ of certain descriptions belonging
to, say, $P(\Z,\Z)$. Comultiplication in $\H_P$
is dual to the composition of descriptions:
$$
\Delta ([p]):=\sum_{q,r|q\circ r=p}[r]\otimes [q].
\eqno(3.1)
$$
(Recall that the composition of descriptions is associative).
 
\smallskip

In order to ensure finiteness of the right hand side of (3.1)
and to produce a Hopf filtration, we must postulate in addition 
existence of a ``size function'' on descriptions. The simplest
properties of such a function $p\mapsto |p|$ that will serve our goal,
are finiteness of the set of descriptions of bounded size,
and additivity
$$
|q\circ r|= |q|+ |r|.
$$
For a concrete example, see [Ma3], subsection 3.4. 

\medskip

{\it (d) Characters, corresponding to the halting problem.}
Finally, assume that we have constructed $\H_P$
and $\Phi (k,f;z)$ as above. Then the character
$\varphi_k :\, \H_P\to \Cal{A}$ (cf. [Ma3], 4.4) corresponding to the halting
problem at a point $k\in \Z$ for the  partial recursive
function computable with the help of a description
$p\in P(\Z,\Z)$, is defined as follows:
$$
\varphi_k ([p]):= \Phi (k,f;z)\in \Cal{A}.
\eqno(3.2)
$$
As soon as this definition is adopted, the machinery and philosophy of
Hopf renormalization and Birkhoff decomposition ([Ma3], Theorem 4.4.1)
becomes applicable to the classical halting problem.
 
 \smallskip
 
Perhaps, it will be even more relevant for quantum computation
schemes based upon infinite--dimensional Hilbert spaces.

\medskip

{\bf 3.2. The simplest construction.} Let $f:\Z\to\Z$ be
a partial recursive function. Consider
its extension $\bar{f}:\,\N\to\N$ defined as follows:
$\bar{f}(x)=f(x)$ if $x\in D(f)$ and $f(x)=0$ otherwise.

\smallskip

Put
$$
\Psi (k,f;z):= \sum_{n=0}^{\infty} \frac{z^n}{(1+n\bar{f}(k))^2}.
\eqno(3.3)
$$
\medskip

{\bf 3.2.1. Proposition.} {\it (i) If $k\notin D(f)$, then
$$
\Psi (k,f;z)=\frac{1}{1-z}.
\eqno(3.4)
$$

(ii) If  $k\in D(f)$, then $\Psi (x,\sigma ;z)$ is the Taylor 
series of a function analytic at $|z|<1$
and continuous at the boundary $|z|=1$. The value
$\bar{f}(k)=f(k)$ can be uniquely reconstructed from
$\Psi$, for example}
$$
f(k)=\sqrt{\left. \frac{dz}{d\Psi}\right|_{z=0}}-1.
\eqno(3.5)
$$

\medskip

The proof is obvious.
\smallskip

Actually, formula (3.3) can be seen in its natural
context if one invokes the general prescription of
reducing any function to a permutation, borrowed from
the theory of quantum computation.
\smallskip

I will briefly recall this prescription now following [Ma3], subsections 3.6 -- 3.8.
\medskip
{\bf 3.3. Reduction of the general halting problem to
the recognition of fixed points of permutations.}
Start with  a partial recursive function $f:\,X\to X$
where $X$ is an infinite constructive world. Extend $X$ by
one point, i.~e. form $X\coprod \{*_X\}$. Choose a total recursive structure
of an additive group without torsion on 
$X\coprod \{*_X\}$ with zero $*_X$.
Extend  ${f}$ to the everywhere defined (but generally uncomputable)
function  $g:\,X\coprod \{*_X\}\to X\coprod \{*_X\}$, by
$$
g(y):= *_X\ \r{if}\ y\notin D({f}).
$$
Define the map
$$
\tau_f:\,(X\coprod \{*_X\})^2 \to (X\coprod \{*_X\})^2
$$
by
$$
\tau_f(x,y):= (x+g(y),y).
\eqno(3.6)
$$
Clearly, it is a permutation. Since $(X\coprod\{*_X\}, +)$
has no torsion, the only finite orbits
of $\tau_f^{\bold{Z}}$ are fixed points. 
\smallskip
Moreover, the restriction of
$\tau_f$ upon the recursive enumerable subset
$$
D(\sigma_f):= (X\coprod \{*_X\})\times D(f)
$$
of the constructive world  $Y:=(X\coprod \{*_X\})^2$
induces a partial recursive permutation $\sigma_f$ 
of this subset. 
\smallskip

Since $g(y)$ never takes the zero value $*_X$ on
$y\in D(f)$, but always is zero outside it,
the complement to $D(\sigma_f)$ in $Y$ consists
entirely of fixed points of $\tau_f$. 

\smallskip

Thus, the halting problem for $f$
reduces to the  fixed point recognition for $\tau_f$. 

\medskip
{\bf 3.4. Permutations with bounded shift.}  The formula (3.3)
can be generalized as follows. 
\medskip

{\bf 3.4.1. Definition.}  {\it Let $\sigma$ be a permutation of $\Z$, $k\in \Z$.
We say that $\sigma$ has a bounded shift at $k$ if there exist
constants $a,b,c$ (depending on $\sigma$ and $k$) such that for all $n\in \bold{Z}$,
$$
c\cdot |n+a|\le \sigma^n(k)\le c\cdot |n+b|.
\eqno(3.7)
$$
}
\medskip

{\bf 3.4.2. Lemma.} {\it  If $\sigma$ has bounded shift at $k$, then the 
$\sigma^{\bold{Z}}$--orbit
of $k$ is infinite, and for any $m\ne 0$ and any point of this orbit $l$,
$\sigma^m$ has bounded shift at $l$.}

\smallskip

{\bf Proof.} Let $l=\sigma^d(k)$. From (3.3) we get
$$
c\cdot |mn+d+a|\le  \sigma^{mn}(l)=\sigma^{mn+d}(k)\le c\cdot |mn+d+b|
$$
that is
$$
c|m|\cdot \left|n+\frac{d+a}{m}\right|\le (\sigma^m)^n(l)\le c|m|\cdot \left|n+\frac{d+b}{m}\right|.
\eqno(3.8)
$$
This inequality has the same form as (3.7), with different constants.

\medskip

{\bf 3.5. Proposition.} {\it Let $\sigma$ be a permutation of $\Z$, $k\in \Z$.
Put
$$
\Psi (k,\sigma ;z):=\sum_{n=1}^{\infty} \frac{z^n}{(\sigma^n(k))^2}.
\eqno(3.9)
$$
Then we have:
\smallskip
(i) If $\sigma^{\bold{Z}}$--orbit of $k$ is finite, then
$\Phi(\sigma ,x;z)$ is a rational function in $z$ whose all poles are of the first order and lie at roots of unity.

\smallskip
(ii) If this orbit is infinite, and $\sigma$ has bounded shift at (any point of) this orbit,
then $\Phi (\sigma ,k;z)$ is the Taylor series of a function analytic at $|z|<1$
and continuous at the boundary $|z|=1$.}

\smallskip

{\bf Proof.}   If $\sigma^{\bold{Z}}$--orbit of $k$ is finite, then
(3.5) is a finite sum of several geometric progressions each of each sums
to a rational function of the type
$$
\r{const}\cdot\frac{z^k}{1-z^l}.
$$
Otherwise, because of (3.7) we get a series absolutely converging
for $|z|\le 1.$ This proves our statement.

\medskip

{\bf 3.6. The Kolmogorov order.} Many interesting $\sigma$,
such as total recursive permutations, are not
permutations of bounded shift.  
 To cope with this situation, we will (uncomputably)
reorder $\Z$, and show, that after this reordering, {\it all} partial
recursive functions and permutations corresponding to them
will satisfy a version of bounded shift property, allowing
one to construct a modification of $\Psi (z)$.

\smallskip

Slightly more generally, let  $X$ be an infinite 
constructive world. Consider an optimal enumeration
$u:\,\Z\to X$ in the sense of Kolmogorov or Schnorr
(see 1.4.4 above). This means that  $u$ is total recursive, surjective,
and the function $C_u:\,X\to \Z$,
$$
C_u(x) := \r{min}\,\{k\,|\,u(k)=x\}
$$
is (a representative of) Kolmogorov complexity of constructive
objects of type $X$.
\smallskip

Now, define the Kolmogorov total order on $X$ associated to $u$ by
$$
x<y\ \Leftrightarrow\ C_u(x)<C_u(y) 
$$
and denote by $\bold{K}=\bold{K}_u:\, X\to \Z$
the function
$$
\bold{K}(x):=1+\r{card}\,\{y\,|\,C_u(y)<C_u(x)\}.
$$
Clearly, $\bold{K}$ is a bijection. If we arrange $X$ in the order of growing Kolmogorov complexity, $\bold{K}(x)$ is precisely the number
of $x$ in this order.

\smallskip

It is convenient also to introduce a Kolmogorov order on $\Z$.
We will denote the respective numbering by the same letter $\bold{K}$.
This should not lead to a confusion.

\smallskip

It is straightforward to check that for some constant $c_0>0$ and all
$x\in X$, we have
$$
c_0\,C_u(x)\le \bold{K}(x)\le C_u(x).
\eqno(3.10)
$$

Let now $\sigma :\,X\to X$ be a partial recursive map,
such that $\sigma$ maps $D(\sigma )$ to $D(\sigma )$
and induces a permutation of this set. Put
$$
\sigma_{\bold{K}} := \bold{K}\circ \sigma\circ \bold{K}^{-1}
$$
and consider this as a permutation of  the subset
$$
D(\sigma_{\bold{K}}):=\bold{K}(D(\sigma )) \subset \Z
$$
consisting of numbers of elements of $D(\sigma )$ in 
the Kolmogorov order. We have then the following modified
version of (3.7):

\medskip

{\bf  3.6.1. Proposition.} {\it Let $x\in D(\sigma)$. If the orbit $\sigma^{\bold{Z}}(x)$
is infinite, then there exist such constants $c_1,c_2 >0$ that
for $k:=\bold{K}(x)$ and all $n\in \bold{Z}$ we have}
$$
c_1\cdot \bold{K}(n)\le \sigma_{\bold{K}}^n(k)\le c_2\cdot \bold{K}(n).
\eqno(3.11)
$$

\smallskip

{\bf Proof.} Let $k=\bold{K}(x),\,x\in X$. We have for $n> 0$:
$$
\sigma_{\bold{K}}^n (k)=\bold{K}(\sigma^n(x))\le c\cdot \bold{K}(n)
\eqno(3.12)
$$
for any fixed Kolmogorov complexity order on $\Z$
(which we denote by the same letter $\bold{K}$ in order to
simplify notation). In fact, if we replace $\bold{K}$ in (3.12)
by the appropriate complexity $C$, this will follow from (2.1), since
$n\mapsto \sigma^n(x)$ is an everywhere defined
morphism $\Z\to X$ in $\c$. 
It remains to invoke (3.10). 

\smallskip

Furthermore, let $Y:=\{\sigma^n(x)\,|\,n\in \Z\}$.
This is a recursively enumerable subset of $X$, and the partial function
$\lambda :\,X\to \Z$ with definition domain $Y$ 
$$
\lambda (y) =n,\ \r{if}\ y=\sigma^n(x)
$$
is partial recursive. Hence again in view of (3.6) and (2.1),
$$
\bold{K}(n)=\bold{K} (\lambda (y))\le c^{\prime}\cdot\bold{K}(y)=
c^{\prime}\cdot\bold{K}(\sigma^n(x)).
\eqno(3.13)
$$
Combining (3.12) and (3.13), we get (3.11) for $n\ge 0.$
Applying the same reasoning to $\sigma^{-1}$ in place of $\sigma$,
we obtain (3.11) for negative $n$. 

\medskip

{\bf 3.7. Proposition.} {\it With the same notations as in Proposition 3.6.1,
put
$$
\Phi (k, \sigma ;z):=\frac{1}{k^2}+\sum_{n=1}^{\infty} \frac{z^{\bold{K}(n)}}{(\sigma_{\bold{K}}^n(k))^2}.
\eqno(3.14)
$$
Then we have:
\smallskip
(i) If $\sigma^{\bold{Z}}$--orbit of $x$ is finite, then
$\Phi(x, \sigma ;z)$ is a rational function in $z$ whose all poles are of the first order and lie at roots of unity.

\smallskip
(ii) If this orbit is infinite,
then $\Phi (x,\sigma ;z)$ is the Taylor series of a function analytic at $|z|<1$
and continuous at the boundary $|z|=1$.}

\medskip

{\bf 3.8. Remarks.} (a) In the proofs of Propositions 3.6.1 and 3.5,
we actually used only the fact that $\sigma$ restricted to the particular
orbit  $Y:=\{\sigma^n(x)\,|\,n\in \Z\}$ is recursive.
thus justifying our choice of $\Cal{A}_-$ in 3.1 (b) above.

\smallskip

(b) Although Kolmogorov's order is as uncomputable as
Kolmogorov's complexity, there are serious arguments
for studying constructions, explicitly involving it, such as
our renormalization characters.

\smallskip

One can argue that all cognitive activity of our
civilization, based upon symbolic (in particular, mathematical)
representations of reality, deals actually with the {\it initial
Kolmogorov segments}
of potentially infinite linguistic constructions,
{\it always} replacing  vast volumes of data by their
compressed descriptions. This is especially visible
in the outputs of the modern genome projects.

\smallskip

In this sense, such linguistic cognitive activity
can be metaphorically compared to a gigantic precomputation 
process, shellsorting infinite worlds of expressions
in their Kolmogorov order.

\bigskip

\centerline{\bf References}

\medskip

[BaSt] J.~Baez, M.~Stay. {\it Physics, topology, logic and computation: a Rosetta stone.} Preprint arxiv:0903.0340

\smallskip

[CalSt1] Ch.~Calude, M.~Stay. {\it Most programs stop quickly or never halt.}
Adv. in Appl. Math., 40 (2008), 295--308
\smallskip

[CalSt2] Ch.~Calude, M.~Stay. {\it Natural halting probabilities,
partial randomness, and zeta functions.} Information
and Computation, 204 (2006), 1718--1739.

\smallskip

[E-FMan] K.~Ebrahimi--Fard and D.~Manchon. {\it The combinatorics of
Bogolyubov's recursion in renormalization.} math-ph/0710.3675

\smallskip

[GaLe] P.~G\'acs, A.~Levin. {\it Causal Nets or What Is a Deterministic Computation?}
Int.~Journ.~Theor.~Phys., vol. 21, No. 12 (1982), 961--971.

\smallskip

[Gr] J.~Grass. {\it Reasoning about Computational Resource Allocation.
An introduction to anytime algorithms.} Posted on the Crossroads website.

\smallskip

[GrZi] J.~Grass, S.~Zilberstein. {\it Programming with anytime algorithms.}
In: Proc. of the IJCAI--95 Workshop on Anytime Algorithms and 
Deliberation Schedulyng. Montreal, 1995.

\smallskip

[He] A.~Heller. {\it An existence theorem for recursive categories.}
Journ. of Symb. Logic, vol. 55, No 3 (1990), 1252--1268.

\smallskip

[Le] L.~Levin. {\it Various measures of complexity for finite objects
(axiomatic description).} Soviet Math. Dokl., vol 17, No. 2 (1976), 522 --526.

\smallskip

[LiVi] Ming Li, P.~Vit\'anyi. {\it An introduction to Kolmogorov complexity
and its applications.} Springer, 1993.

\smallskip

[Ma1]  Yu.~Manin.  {\it A Course in Mathematical Logic. } Springer Verlag, 1977. XIII+286 pp. (The second, expanded Edition to be published in 2009). 

\smallskip

[Ma2]  Yu.~Manin. {\it Classical computing, quantum computing,
and Shor's factoring algorithm.}  S\'eminaire Bourbaki, no. 862 (June 1999),
Ast\'erisque, vol 266, 2000, 375--404.
 quant-ph/9903008.
 
 \smallskip
 
 [Ma3] Yu.~Manin. {\it Renormalization and computation I. Motivation and background.}
Preprint arxiv:0904.492

\smallskip

[NaWe] A. Nabutovsky, S.~Weinberger. {\it The fractal nature of
Riemm/Diff I.} Geometriae Dedicata, 101 (2003), 145--250.

\smallskip

[Ro] H.~Rogers. {\it G\"odel numberings of partial recursive functions.}
Journ. Symb. Logic, 23 (1958), 331--341.

\smallskip

[RuZi]  S.~J.~Russell, S.~Zilberstein. {\it Composing real--time systems.} In:
Proc. of the XIIth International Joint Conference on Artificial Intelligence,
(1991), Sydney, pp. 212--217.

\smallskip

[Sch] C.~P.~Schnorr. {\it Optimal enumerations and optimal G\"odel
numberings.} Math. Systems Theory, vol. 8, No. 2 (1974), 182--191.

\smallskip

[So]  R.~I.~Soare. {\it Computability theory and differential geometry.}
Bull. of Symb. Logic, vol. 10, Nr 4 (2004), 457--486.

\smallskip

[Ya] N.~S.~Yanofsky. {\it Towards a definition of an algorithm.} math.LO/0602053

\enddocument

\centerline{\bf Glossary}

\medskip

{\bf Bit size.} My bit size of $n$ is Schnorr's $||n||$.

\medskip

{\bf Effective.} Synonymous with total recursive, computable 

{\bf Enumeration.} In [Sch], p. 182: optimal enumerations of
recursively enumerable sets. (it seems that {\it enumerations}
refer to sets, whereas {\it numberings} to functions. 

\smallskip

In [Sch], p. 187: optimal enumeration of $\N$ is the same as optimal
G\"odel numbering  of functions of arity 0.

\medskip

{\bf  Numbering.} My current use:  
bijection between $\Z$ and elements of a constructive world.
In the book, I allow decidable sets. Rogers [Ro] as well allows decidable sets.

\smallskip

[Ro] works with $\N$ and numbers descriptions of p.~r. functions (of 1 variable0. 
More precisely: (unqualified) {\it numbering} is a versal family
of p.r. functions, containing {\it all} such functions.
It is {\it semi--effective}, if it is induced by a p.r. function of 2 variables
over the initial base.
It is {\it fully effective} if

In [Sch]: G\"odel numbering of partial recursive
functions of arity $k$ is a  "strongly versal" family of such functions.
Strong versality means that any family is induced from it by a total
recursive morphism of bases.

\smallskip

I will use the word {\it numberings} quite differently:  as computable 
together with their inverses
{\it bijective} maps, or as in my book bijections with decidable subsets of $\N$.
\smallskip
Schnorr's numberings of functions of arity $k$ are my versal
families.

\medskip

{\bf Optimal G\"odel numbering.} It is such a strongly versal family, that
totally recursive base change to any other family can be taken linearly bounded.

\smallskip

{\bf NB} Main Schnorr's results:  {\it a) strongly versal families exist.
b) They are related to each other by linearly bounded total 
recursive permutations. }

\smallskip

Kolmogorov's optimality is {\it strictly weaker}! This is Schnorr's 
Theorem 10.

\smallskip

This refers also to arity 0, and probably to any arity. Schnorr
for some reason treats separately cases of arity 1 and 0 only.

\medskip

{\bf Program.} For me a program has variable inputs and outputs,
so it is something like a Turing machine. In many other papers,
programs is what I call {\it (an instance of) an input}. So when they
say {\it a program}, they implicitly assume a
fixed Turing machine etc.

\medskip

{\bf Program complexity.} My complexity $C$.

\medskip

{\bf Recursive.} [Sch] means by this {\it total recursive.}

\bigskip

\centerline{\bf On natural numbers objects}

\medskip

$\bullet$ One could start with $\Cal{N}$, and then postulate $\times$ and $\coprod$
and that everything is isomorphic to $\Cal{N}$. One then can even postulate
that $\Cal{N}=s(\Cal{N})\coprod \bold{1}$ and $s:\,\Cal{N}\to s(\Cal{N})$ is an isomorphism.

\medskip

$\bullet$ Something like $\Cal{N}$ is needed in order to make
invariant sense of ``linearly bounded'' morphisms. Since my numberings
are defined only up to total recursive renumberings, we have to define
intrinsically the conjugacy class of subgroups of linearly bounded growth.

\enddocument